\documentclass[12pt]{article}
\usepackage{amssymb}
\usepackage{mathrsfs}%Óиü¶àµÄÊýѧ·ûºÅ£¬±ÈÈ绨Ìå$\mathscr P$
\usepackage{pst-poly}     % From pstricks/contrib/pst-poly
\usepackage{cite}

\newtheorem{thm}{Theorem}[section]
\newtheorem{lem}[thm]{Lemma}
\newtheorem{Def}[thm]{Definition}
\newtheorem{cor}[thm]{Corollary}

\newenvironment{pf}[1][Proof]{\noindent\textbf{#1.} }{\hfill\rule{1mm}{2mm}}

\makeatletter \@addtoreset{equation}{section} \makeatother

\begin{document}
\title{\bf Fault Diagnosability of Arrangement Graphs \\
\baselineskip 12pt \vskip0.05cm
\author{Shuming Zhou$^1$ \quad Jun-Ming Xu$^2$\\
\small $^1$College of Mathematics and Computer Science, Fujian
Normal
University, \\
\small Fuzhou, Fujian, P.R.China 350007  \\
\small $^2$Department of Mathematics, University of Science and Technology of China, \\
\small Hefei, Anhui, 230026, P.R. China\\
 }}\date{}
 \maketitle

\footnotetext{\small *This work was partly supported by National
Natural Science foundation of China (No. 61072080, 11071233) and the
Key Project of Fujian Provincial Universities services to the
western coast of
the straits-Information Technology Research Based on Mathematics.\\
E-mail address: zhoushuming@fjnu.edu.cn}

{\bf\noindent Abstract:}\ The growing size of the multiprocessor
system increases its vulnerability to component failures. It is
crucial to locate and to replace the faulty processors to maintain
a system's high reliability. The fault diagnosis is the process of
identifying faulty processors in a system through testing. This
paper shows that the largest connected component of the survival
graph contains almost all remaining vertices in the $(n,
k)$-arrangement graph $A_{n,k}$ when the number of moved faulty
vertices is up to twice or three times the traditional
connectivity. Based on this fault resiliency, we establishes that
the conditional diagnosability of $A_{n,k}$ under the comparison
model. We prove that for $k\geq 4$, $n\geq k+2$, the conditional
diagnosability of $A_{n,k}$ is $(3k-2)(n-k)-3$; the conditional
diagnosability of $A_{n,n-1}$ is $3n-7$ for $n\geq 5$.

\vskip6pt

{\noindent{\bf Keywords:} Fault tolerance; comparison diagnosis;
diagnosability; $(n, k)$-arrangement graph.

\section{Introduction }

Distributed processor architectures offer the potential advantage
of high speed, provided that they are highly fault-tolerant and
reliable, and have good communication between remote processors.
An important component of such a distributed system is its network
topology, which defines the inter-processor communication
architecture. Fault-tolerance is especially important for
interconnection networks, since computers may fail, creating
faults in the network. To be reliable, the rest of the network
should stay connected. Obviously, this can only be guaranteed if
the number of faults is smaller than the minimum degree in the
network. When the number of faults is larger than the minimum
degree, some extensions of connectivity are necessary, since the
graph may become disconnected. Some generalizations of
connectivity were introduced and examined for various classes of
graphs in~\cite{clp01}, including super connectedness and tightly
super connectedness, where only one singleton can appear in the
remaining network, and restricted connectivity and super
connectivity, where a remaining component must have a certain
minimum size. As we increase the number of faults in the graph, it
is desirable that the largest component of the surviving network
stays connected, with a few processors separated from the rest,
since then the network will continue to be able to function. Many
interconnection networks have been examined in this aspect, when
the number of faults is roughly twice the minimum degree,
see~\cite{cl07a, kllht09}. One can even go further and ask what
happens when more vertices are deleted. This has been examined for
the hypercube in~\cite{yecml04, yem04, yem06} and for certain
Cayley graphs generated by transpositions in~\cite{cl07b}, and it
has been shown that the surviving network has a large component
containing almost all vertices.

The process of identifying faulty processors in a system by
analyzing the outcomes of available inter-processor tests is called
system-level diagnosis. In 1967, Preparata, Metze, and
Chien~\cite{pmc67} established a foundation of system diagnosis and
an original diagnostic model, called the PMC model. Its target is to
identify the exact set of all faulty vertices before their repair or
replacement. All tests are performed between two adjacent
processors, and it was assumed that a test result is reliable
(respectively, unreliable) if the processor that initiates the test
is fault-free (respectively, faulty). The comparison-based diagnosis
models, first proposed by Malek~\cite{m80} and Chwa and
Hakimi~\cite{ch81}, have been considered to be a practical approach
for fault diagnosis in the multiprocessor systems. In these models,
the same job is assigned to a pair of processors in the system and
their outputs are compared by a central observer. This central
observer performs diagnosis using the outcomes of these comparisons.
Maeng and Malek~\cite{mm81} extended Malek's comparison approach to
allow the comparisons carried out by the processors themselves.
Sengupta and Dahbura\cite{s92} developed this comparison approach
such that the comparisons have no central unit involved.

Lin et al.~\cite{lthcl08} introduced the conditional diagnosis
under the comparison model. By evaluating the size of connected
components, they obtained that the conditional diagnosability of
the star graph $S_n$ under the comparison model is $3n-7$, which
is about three times larger that the classical diagnosability of
star graphs. In the same method, Hsu et al.~\cite{hcsht09} have
recently proved that the conditional diagnosability of the
hypercube $Q_n$ is $3n-5$. This idea was attributed to Lai et
al.~\cite{ltch05} who are the first to use a conditional diagnosis
strategy. They obtained that the conditional diagnosability of the
hypercube $Q_n$ is $4n-7$ under the PMC model. Furthermore, Hsu et
al.~\cite{hcsht09} exposed the difference between these two
conditional diagnosis models.

The arrangement graph, proposed as a generalization of the star
graph in an attempt to solve the scalability problem of the star
graph topology, while preserving its attractive features, has been
extensively studied~\cite{benm98, cjs01, cjt01, cc98, dt91, dt92,
dt93, hlth04, hch99, lc01, tmh08}. Based on the fault tolerance of
the arrangement graph,
%Based on the fault tolerance of the $(n,k)$-arrangement graph
%$A_{n,k}$ ($n\geq 4$, $k\geq 3$, $n-k\geq 1$),
in this paper, we establish its conditional diagnosability under the
comparison diagnosis model. The rest of this paper is organized as
follows. Section 2 introduces some definitions, notations and the
structure of the arrangement graph. Section 3 is devoted to the
fault resiliency of $A_{n,k}$, and Section 4 concentrates on the
conditional diagnosability of the arrangement graph. Section 5
concludes the paper.

%\section{Preliminaries}
\section{Arrangement graphs}

For notation and terminology not defined here we follow
\cite{x01}. Specifically, we use a graph $G=G(V,E)$ to represent
an interconnection network, where a vertex $u\in V$ represents a
processor and an edge $(u,v)\in E$ represents a link between
vertices $u$ and $v$. If at least one end of an edge is faulty,
the edge is said to be faulty; otherwise, the edge is said to be
fault-free. Let $S$ be a subset of $V(G)$. The subgraph of $G$
induced by $S$, denoted by $G[S]$, is the graph with the
vertex-set $S$ and the edge-set $\{(u,v)|\ (u,v)\in E(G), u,v\in
S\}$. For a vertex $u$ in $G$, $N(u)$ denotes the set of all
neighbors of $u$, i.e., $N(u)=\{v|\ (u,v)\in E\}$. Let $S$ be a
subgraph of $G$ or a subset of $V(G)$, and let $N(S)=\bigcup_{u\in
S}(u)\setminus S$.
% , $N[F]=N(F)\cup F$. For brevity, $N[u]=N(u)\cup \{u\}$, $N(\{u,v\})$
%and $N[\{u,v\}]$ are written as $N(u,v)$ and $N[u,v]$, respectively.
%The symmetric difference of two sets $F_1$ and $F_2$ is defined as
%the set $F_1\bigtriangleup F_2=(F_1-F_2)\cup (F_2-F_1)$.
We use $K_n$ to denote the complete graph of order $n$, and $d(u,v)$
to denote the distance between $u$ and $v$, the length of a shortest
path between $u$ and $v$ in $G$. The diameter of $G$ is defined as
the maximum distance between any two vertices in $G$.

For any subset $F\subset V$, the notation $G-F$ denotes a graph
obtained by removing all vertices in $F$ from $G$ and deleting those
edges with at least one end-vertex in $F$, simultaneously. If $G-F$
is disconnected, $F$ is called a {\it separating set}. A separating
set $F$ is called a {\it $k$-separating set} if $|F|=k$. The maximal
connected subgraphs of $G-F$ are called {\it components}. The {\it
connectivity} $\kappa (G)$ of $G$ is defined as the minimum $k$ for
which $G$ has a $k$-separating set; otherwise $\kappa (G)$ is
defined $n-1$ if $G=K_n$. A graph $G$ is called to be {\it
$k$-connected} if $\kappa (G)\geq k$. A $k$-separating set is called
to be {\it minimum} if $k=\kappa (G)$.

The interconnection network has been an important research area for
parallel and distributed computer systems. Network reliability is
one of the major factors in designing the topology of an
interconnection network. The well-known hypercube is
% and its variants are
the first major class of interconnection networks.

As another topology of an interconnection network, Akers and
Krishnamurthy~\cite{ak89} proposed the star graph $S_n$, which has
smaller degree, diameter, and average distance than the comparable
hypercube, while reserving symmetry properties and desirable
fault-tolerant characteristics. As a result, the star graph has
been recognized as an alternative to the hypercube. However, the
star graph is less flexible in adjusting its sizes. With the
restriction on the number of vertices, there is a large gap
between $n!$ and $(n+1)!$ for expanding an $S_n$ to $S_{n+1}$. To
relax the restriction of the numbers of vertices $n!$ in $S_n$,
The arrangement graph was proposed by Day and Tripathi~\cite{dt92}
as a generalization of the star graph $S_n$. It is more flexible
in its size than $S_n$.

\begin{Def}
Given two positive integers $n$ and $k$ with $n>k$, let $\langle
n\rangle$ denote the set $\{1,2,\ldots, n\}$, and let $P_{n,k}$ be a
set of arrangements of $k$ elements in $\langle n\rangle$. The {\it
$(n,k)$-arrangement graph}, denoted by $A_{n,k}$, has vertex-set
$V(A_{n,k})=P_{n,k}$ and edge-set $E(A_{n,k})=\{(p,q)\ |$\ $p$ and
$q$ differ in exactly one position $\}$.
\end{Def}

The graph shown in Figure~\ref{f1} is a $(4,2)$-arrangement graph
$A_{4,2}$.

\vskip6pt

\begin{figure}[h]
\begin{pspicture}(-3.25,0)(7.5,6.5)
\psset{radius=.1} \Cnode(.8,1){43}\rput(0.7,.7){\scriptsize 4\,3}
\Cnode(2.5,1){23}\rput(2.5,.7){\scriptsize 2\,3}
\Cnode(5.5,1){21}\rput(5.5,.7){\scriptsize 2\,1}
\Cnode(7.2,1){41}\rput(7.3,.7){\scriptsize 4\,1}
\Cnode(1.75,2.5){13}\rput(1.4,2.5){\scriptsize 1\,3}
\Cnode(4,2){24}\rput(4,1.7){\scriptsize 2\,4}
\Cnode(6.25,2.5){31}\rput(6.6,2.5){\scriptsize 3\,1}
\Cnode(3.3,3.5){14}\rput(3,3.7){\scriptsize 1\,4}
\Cnode(4.7,3.5){34}\rput(5.1,3.6){\scriptsize 3\,4}
\Cnode(4.7,4.7){32}\rput(5.1,4.7){\scriptsize 3\,2}
\Cnode(3.3,4.7){12}\rput(2.9,4.7){\scriptsize 1\,2}
\Cnode(4,6){42}\rput(4,6.3){\scriptsize 4\,2}

\ncline[linecolor=red,linewidth=1.5pt]{43}{23}\ncline{21}{23}\ncline[linecolor=red,linewidth=1.5pt]{41}{21}\ncline[linecolor=red,linewidth=1.5pt]{43}{13}
\ncline[linecolor=red,linewidth=1.5pt]{13}{23}\ncline{24}{23}\ncline{24}{21}\ncline[linecolor=red,linewidth=1.5pt]{21}{31}
\ncline{13}{12}\ncline{13}{14}\ncline{42}{12}\ncline[linecolor=red,linewidth=1.5pt]{14}{24}
\ncline[linecolor=red,linewidth=1.5pt]{24}{34}\ncline{34}{31}\ncline[linecolor=red,linewidth=1.5pt]{14}{34}\ncline{14}{12}
\ncline{32}{34}\ncline{32}{31}\ncline[linecolor=red,linewidth=1.5pt]{12}{32}\ncline[linecolor=red,linewidth=1.5pt]{12}{42}
\ncline[linecolor=red,linewidth=1.5pt]{42}{32}\ncline[linecolor=red,linewidth=1.5pt]{31}{41}
\nccurve[angleA=-160,angleB=100]{42}{43}
\nccurve[angleA=-20,angleB=80]{42}{41}
\nccurve[angleA=-30,angleB=-150]{43}{41}
\end{pspicture}
\caption{\label{f1}\footnotesize{The structure of $A_{4,2}$}}

\end{figure}
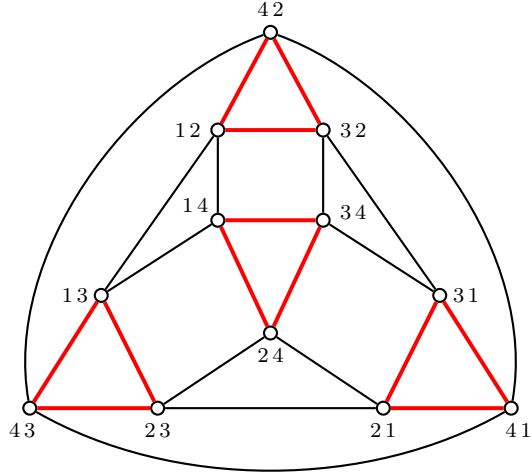

Clearly, $A_{n,k}$ is a $k(n-k)$-regular graph with
$\frac{n!}{(n-k)!}$ vertices. It was showed by Day and
Tripathi~\cite{dt92} that $A_{n,k}$ is vertex-symmetric and
edge-symmetric and has the diameter of $\lfloor
\frac{3k}{2}\rfloor$. Day and Tripathi~\cite{dt91} showed
%that there exist $k(n-k)$ vertex-disjoint paths between any two vertices in the arrangement
%graph, that is,
the connectivity $\kappa(A_{n,k})=k(n-k)$.

Moreover, $A_{ n,1}$ is isomorphic to the complete graph $K_n$, and
$A_{ n,n-1}$ is isomorphic to the $n$-dimensional star graph $S_n$.
Chiang and Chen~\cite{cc98} showed that $A_{n,n-2}$ is isomorphic to
the $n$-alternating group graph $AG_n$.
%, and they also obtained the average distance of $A_{n,k}$ and $AG_n$.

For two distinct $i$ and $j$ in $\langle n\rangle$, let
$V_{n,k}^{j:i}$ be the set of all vertices in $A_{n,k}$ with the
$j$th position being $i$, that is,
 $$
 V_{n,k}^{j:i}=\{p\, |\ p=p_1\cdots p_j\cdots p_k\in P_{n,k}\ {\rm and}\ p_j=i \}.
 $$
For a fixed position $j\in \langle n\rangle$, $\{V_{n,k}^{j:i}|\
1\leq i\leq n\}$ forms a partition of $V(A_{n,k})$. Let
$A_{n,k}^{j:i}$ denote the subgraph of $A_{n,k}$ induced by
$V_{n,k}^{j:i}$. Then for each $j\in \langle n\rangle$,
$A_{n,k}^{j:i}$ is isomorphic to $A_{n-1,k-1}$. For example, a
partition of $A_{4,2}$ is shown in Figure~\ref{f1}, where red
triangles are $A_{4,2}^{2:i}$'s with $i\in \langle 4\rangle$,
isomorphic to $A_{3,1}=K_3$.

Thus, $A_{n,k}$ can be recursively constructed from $n$ copies of
$A_{n-1,k-1}$. It is easy to check that each $A_{n,k}^{j:i}$ is a
subgraph of $A_{n,k}$, and we say that $A_{n,k}$ is decomposed into
$n$ subgraphs $A_{n,k}^{j:i}$'s according to the $j$th position. For
simplicity, by the symmetry of $A_{n,k}$ we shall take $j$ as the
last position $k$, and use $A_{n,k}^{i}$ to denote $A_{n,k}^{k:i}$.

Let $E(i,j)$ be the set of edges between $A_{n,k}^{i}$ and
$A_{n,k}^{j}$, that is,
 $$
 E(i,j)=\{(p,q)\in E(A_{n,k})|\ p\in V(A_{n,k}^{i})\ {\rm and}\ q\in (A_{n,k}^{j})\}.
 $$
Clearly, $E(i,j)$ is a perfect matching ( a set of edges in which
any two edges have no common end-vertex) between $A_{n,k}^i$ and
$A_{n,k}^j$, and
 \begin{equation}\label{e2.1}
 |E(i,j)|=\frac{(n-2)!}{(n-k-1)!}.
 \end{equation}

Let $I$ be a subset of $\langle n\rangle$, and let $H$ be a subset
of $V(A_{n,k}^I)$ or a subgraph of $A_{n,k}^I$, where
$A_{n,k}^I=\{A_{n,k}^i:\ i\in I\}$.
% for some $i\in\langle n\rangle$ and .
Use $N^I(H)$ to denote the set of neighbors of $H$ in $A_{n,k}^I$.
Particularly, use $N^{\overline{I}}(H)$ and $N^{I}(H)$ as an
abbreviation of $N^{\langle n\rangle\setminus{I}}(H)$ and
$N^{I}(H)$, respectively, and call vertices in
$N^{\overline{I}}(H)$ and $N^{I}(H)$ the {\it outer neighbors} and
{\it inner neighbors} of $H$, respectively. Obviously, every
vertex $u$ of $A^i_{n,k}$ has $n-k$ outer neighbors, and two
arbitrary outer neighbors of $u$ are distributed in distinct
subgraphs. We write $u$ for $\{u\}$. It follows from the
definitions that, for every $i\in\langle n\rangle$,
 \begin{equation}\label{e2.2}
 \begin{array}{rl}
 &|N^{i}(u)|=(k-1)(n-k)\ {\rm and}\ |N^{\overline{i}}(u)|=n-k,\\
 & N^{\overline{i}}(u)\cap N^{\overline{i}}(v)=\emptyset\ \ {\rm if}\ u, v\in A_{n,k}^i\ {\rm and }\ u\ne
 v,
 \end{array}
 \end{equation}
and for any two distinct vertices $x\in A_{n,k}^i$ and $y\in
A_{n,k}^j$ with $i\ne j$, and $I=\{i,j\}$,
\begin{equation}\label{e2.3}
 |N^{\overline{I}}(x)\cap N^{\overline{I}}(y)|=0\ \ {\rm if}\ x\ {\rm and}\ y\ {\rm are\ not\ adjacent}.
 \end{equation}

We say that one vertex $u$ is adjacent to some subgraph $A_{n,k}^j$
if $u$ has an outer neighbor in $A_{n,k}^j$. Let
 $$
 V_i=\{u_1u_2\cdots
u_{i-1}xu_{i+1}\cdots u_k|\ x\in \langle n\rangle\setminus\{u_1,
u_2, \cdots , u_{i-1}, u_{i+1}, \cdots , u_k\}\}
$$
Then, when $n\geq k+2$, the graph induced by $V_i$ is a complete
graph of order $n-k+1$ and a subgraph of $A_{n,k}^{u_k}$, which
implies that any two adjacent vertices have exactly $(n-k-1)$ common
neighbors. Thus, by the edge-transitiveness of $A_{n,k}$, for any
edge $e$,
 \begin{equation}\label{e2.4}
 |N(e)|=2k(n-k)-(n-k-1)-2=(2k-1)(n-k)-1.
 \end{equation}

In addition, the following property of $A_{n,k}$ is useful, which
can be checked by the definition of $A_{n,k}$. For any two distinct
vertices $u$ and $v$ in $A_{n,k}$,
 \begin{equation}\label{e2.5}
 |N(u)\cap N(v)|=\left\{
 \begin{array}{ll}
 0,\ \ &{\rm if}\ d(u,v)\geq 3;\\
 2,\ \ &{\rm if}\ d(u,v)=2\ {\rm and}\ n\geq k+2;\\
 1,\ \ &{\rm if}\ d(u,v)=2\ {\rm and}\ n=k+1;\\
 n-k-1,\ \ &{\rm if}\ d(u,v)=1.
 \end{array}\right.
 \end{equation}

Other properties of the arrangement graph has received
considerable attention in the literature. First, Day and
Tripathi~\cite{dt93} showed the existence of pancyclicity, that is
$A_{n,k}$ contains cycles of all lengths. Hsieh et
al.~\cite{hch99} investigated the existence of hamiltonian cycle
in $A_{n,k}$ with faulty vertices, Lo and Chen~\cite{lc01} studied
hamiltonian connectedness of $A_{n,k}$ with faulty edges. Hsu et
al.~\cite{hlth04} further obtained an optimal result that the
graph $A_{n,k}$ ($n\geq k+2$) is $(k(n-k)-2)$-hamiltonian and
$(k(n-k)-3)$-hamiltonian connected
%(where, a graph $G$ is $f$-hamiltonian (respectively,
%$f$-hamiltonian connected) if there is a hamiltonian cycle
%(respectively, a hamiltonian path between any two vertices)
in $G-F$ for any $F\subset V(G)\cup E(G)$ with $|F|\leq f$). Teng
et al.~\cite{tmh08} have recently shown that $A_{n,k}$ is
panpositionable hamiltonian and panconnected if $k>1$ and $n\geq
k+2$
%(where, a hamiltonian graph $G$ is panpositionable
%(respectively, panconnected) if for any two different vertices $u$
%and $v$ of $G$, there exists a hamiltonian cycle $C$(respectively,
%path $P$) joining any two different vertices $u$ and $v$ of $G$
% such that the relative distance between $u$ and $v$ on $C$
% (respectively, $P$) is $l$, where, the integer $l$ satisfying
% $d(u,v)\leq l\leq |V(G)|-d(u,v)$ (respectively, $d(u,v)\leq l\leq
% |V(G)|-1$))
.
In addition, Bai et al.~\cite{benm98} proposed a distributed
algorithm with optimal time complexity and without message
redundancy for one-to-all broadcasting in one-port communication
model on the fault-free arrangement graphs, and also developed a
fault tolerant broadcasting algorithm with less than $k(n-k)$
faulty edges. Chen et al.~\cite{cjs01,cjt01} presented efficient
one/all-to-all broadcasting algorithms on the arrangement graphs
by constructing $n-k$ spanning trees, where the height of each
tree is $2n-1$.

%The recursively constructed properties of the arrangement graph will
%be given in the next section.

\section{Fault tolerance of the arrangement graph}

The connectivity $\kappa (G)$ of a graph $G$ is an important
parameter to measure the fault tolerance of the network, while it
has an obvious deficiency in that it tacitly assume that all
elements in any subset of $G$ can potentially fail at the same time.
To compensate for this shortcoming, it would seem natural to
generalize the classical connectivity by introducing some conditions
or restrictions on the separating set $S$ and/or the components of
$G-S$.

Recall the connectivity $\kappa (G)$ of $G$, it is the minimum
number of vertices whose removal results in a disconnected or a
trivial (one vertex) graph. A $k$-regular $k$-connected graph is
%A $k$-regular graph is {\it maximally connected} if it is
%$k$-connected. A regular graph is
{\it super $k$-connected} if any one of its minimum separating sets
is a set of the neighbors of some vertex. If, in addition, the
deletion of a minimum separating set results in a graph with two
components (one of which has only one vertex), then the graph is
{\it tightly super $k$-connected}. For example, the complete
bipartite graph $K_{n,n}$ is $n$-super connected but not tightly
$n$-super connected. The notions of super connectedness and tightly
super connectedness were first introduced in~\cite{bbst81}
and~\cite{clp01}, respectively.

Esfahanian~\cite{e89} first introduced the concepts of the
restricted separating set and the restricted connectivity of a graph
$G$. A set $S$ of vertices is a {\it restricted separating set} if
$G-S$ is disconnected and $N(x)$ is not completely contained in $S$
for any vertex $x$ in $G$. The {\it restricted connectivity} of $G$,
denoted by $\kappa_r(G)$, is the minimum cardinality of a restricted
vertex-cut.

Considering it is not easy to examine whether a separating set is
restricted, Xu et al.~\cite{xww10} formally proposed the super
connectivity, a weaker concept than the restricted connectivity. A
separating set $S$ of $G$ is {\it super} if $G-S$ contains no
isolated vertices. The {\it super connectivity} of $G$, denoted by
$\kappa_s(G)$, is the minimum cardinality of a super separating set.
Clearly, $\kappa(G)\leq\kappa_s(G)\leq\kappa_r(G)$ if $\kappa_r(G)$
exists.

It follows from definitions that the restricted connectivity or
super connectivity can provide a more accurate measurement than the
connectivity for fault tolerance of a large-scale interconnection
network.
%It has been
%shown that if a network is restricted connected, it is more reliable
%and has the lower vertex failure comparing to that has only the
%super connectivity property.

Usually, if the surviving graph $G-S$ contains a large connected
component $C$ when $G-S$ is not connected, the component $C$ may be
used as the functional subsystem, without incurring severe
performance degradation. Thus, in evaluating a distributed system,
it is indispensable to estimate the size of the maximal connected
components of the underlying graph when the structure begins to lose
processors.

Yang et al.~\cite{yecml04, yem04, yem06} proved that the hypercube
$Q_n$ with $f$ faulty processors has a component of size $2^n-f-1$
if $f\leq 2n-3$, and size $2^n-f-2$ if $f\leq 3n-6$. Yang et
al.~\cite{ymtx08, ymlc07} also obtained that a similar result for
the star graph $S_n$. Cheng et al.~\cite{cl02, cl06} gave a more
detail result for $S_n$. The removal of any separating set of at
most $2n-4$ from $S_n$ results in exact two components, one of
them is a single vertex or edge. Cheng and Lipt\'ak~\cite{cl07b}
generalized this result for $S_n$ with linearly many faults. Cheng
et al.~\cite{cls10} presented a similar result for in
2-tree-generated networks with linearly many faults. In this
section, we detail on the fault resilience of the arrangement
graph $A_{n,k}$.

Throughout this paper, the notation $F$ denotes a set of vertices in
$A_{n,k}$. If $F$ is regarded as a set of faulty vertices, then a
subgraph $H$ of $A_{n,k}$ is called {\it fault-free} if $V(H)\cap
F=\emptyset$. Let

\begin{equation}\label{e3.1}
 F_i=A_{n,k}^i\cap F\ \ {\rm and}\ \ f_i=|F_i|\ {\rm for}\ 1\leq i\leq n.
 \end{equation}

We first discuss the tightly super connectedness. Since $A_{n,1}$
is isomorphic to a complete graph $K_n$, it is super connected but
not tightly super connected. When $n=4$, it is easy to check that
$A_{4,2}$ is not tightly super connected since it has a separating
set $F$ with $|F|=4$ such that two components of $A_{4,2}-F$ are
both 4-cycles (see Figure~\ref{f1}). Thus, in the following
discussion, we assume $k\geq 3$.
%$3\leq k\leq n-2$.

\vskip6pt

\begin{thm}\label{thm3.1}
For $k\geq 3$, $A_{n,k}$ is tightly super $k(n-k)$-connected.
\end{thm}

\begin{pf}
Let $F$ be a minimum separating set in $A_{n,k}$. Then, using the
notations defined in (\ref{e3.1}), we have that
 $$
 |F|=\sum^n_{i=1}f_i=\kappa(A_{n,k})=k(n-k).
 $$
By the definition of tightly super connectivity, we need to show
that $A_{n,k}-F$ has exactly two components, one of them is a single
vertex. We gain our ends by proving the following claims.

\vskip6pt

{\bf Claim 3.1.1}\ $f_i\geq (k-1)(n-k)$ for some $i\in \langle
n\rangle$.

{\it Proof:}\ Suppose to the contrary that $f_i<(k-1)(n-k)$ for any
$i\in \langle n\rangle$. Then $A_{n,k}^i-F_i$ is connected since
$A_{n,k}^i$ is $(k-1)(n-k)$-connected. We will deduce a
contradiction by showing that $A_{n,k}-F$ is connected. To this end,
we only need to show that $A_{n,k}^i$ and $A_{n,k}^j$ can be
connected in $A_{n,k}-F$ for any two distinct $i, j\in \langle
n\rangle$.

In fact, by (\ref{e2.1}), when either $k\geq 4$ or $k=3$ and $n\geq
6$, we have
 $$
|E(i,j)|=\frac{(n-2)!}{(n-k-1)!}>k(n-k)=|F|,
 $$
which implies that there exists a fault-free edge $e\in E(i,j)$. It
follows that $A_{n,k}^i$ and $A_{n,k}^j$ can be connected in
$A_{n,k}-F$ by the fault-free edge $e$.

When $k=3$ and $n\in\{4,5\}$, we have
 \begin{equation}\label{e3.2}
 |E(i,j)|=\frac{(n-2)!}{(n-k-1)!}=\left\{\begin{array}{ll}
 2<3=|F|\ &\ {\rm if}\ n=4;\\
 6=|F|\ &\ {\rm if}\ n=5.
 \end{array}\right.
 \end{equation}
Without loss of generality, assume that there are no fault-free
edges in $E(i,j)$ (otherwise $A_{n,3}^i$ and $A_{n,3}^j$ can be
connected in $A_{n,3}-F$ by some fault-free edge in $E(i,j)$). By
(\ref{e3.2}), there exist a fault-free edge $e_1$ in $E(i,x)$ and a
fault-free edge $e_2$ in $E(x,j)$ for any $x\notin\{i,j\}$. Thus,
$A_{n,3}^i$ and $A_{n,3}^j$ can be connected in $A_{n,3}-F$ by
$A_{n,3}^x$ and the fault-free edges $e_1$ and
$e_2$.\hfill\rule{1mm}{2mm}

 \vskip6pt

{\bf Claim 3.1.2}\ If there is some $i\in \langle n\rangle$ such
that $|F-F_i|<(k-1)(n-k)$, then $A_{n,k}-(V(A_{n,k}^i)\cup (F-F_i))$
is connected.
%In particular, if $f_i=(k-1)(n-k)$ then $A_{n,k}^i-F_i$ is disconnected.

{\it Proof:}\ By the hypothesis, for any $j\in \langle n\rangle$
with $j\ne i$, we have $f_j<(k-1)(n-k)$, which implies that
$A_{n,k}^j$ is connected since $A_{n,k}^j$ is
$(k-1)(n-k)$-connected. Since for any two distinct $j,t\in\langle
n\rangle\setminus\{i\}$,
 $$
|E(j,t)|=\frac{(n-2)!}{(n-k-1)!}>(k-1)(n-k)>|F-F_i|,
 $$
there exists a fault-free edge $e$ in $E(j,t)$. Thus $A_{n,k}^j$ and
$A_{n,k}^t$ can be connected in $A_{n,k}-F$ by the fault-free edge
$e$. By the arbitrariness of $j$ and $t$, $A_{n,k}-(V(A_{n,k}^i)\cup
(F-F_i))$ is connected. \hfill\rule{1mm}{2mm}

\vskip6pt

{\bf Claim 3.1.3}\ $f_i\leq (k-1)(n-k)$ for any $i\in \langle
n\rangle$.

{\it Proof:}\ If there is some $i\in \langle n\rangle$ such that
$f_i>(k-1)(n-k)$, then
 $$
 |F-F_i|<k(n-k)-(k-1)(n-k)=n-k<(k-1)(n-k).
 $$
By Claim 3.1.2, $A_{n,k}-(V(A_{n,k}^i)\cup (F-F_i))$ is connected.
Since every vertex in $A_{n,k}^i-F_i$ has exactly $n-k$ outer
neighbors in $A_{n,k}-A_{n,k}^i$ and $|F-F_i|<n-k$, and at least one
of the $n-k$ outer neighbors is fault-free, $A_{n,k}-F$ is still
connected, a contradiction. \hfill\rule{1mm}{2mm}

\vskip6pt

We now show our theorem. By Claim 3.1.1 and Claim 3.1.3, there
exists some $i\in \langle n\rangle$ such that $f_i=(k-1)(n-k)$.
Thus, for $k\geq 3$,
 $$
|F-F_i|=k(n-k)-(k-1)(n-k)=n-k<(k-1)(n-k).
 $$
By Claim 3.1.2, $A_{n,k}-(A_{n,k}^i\cup (F-F_i))$ is connected,
which implies $A_{n,k}-A_{n,k}^i$ is $(n-k+1)$-connected.

Suppose that  $A_{n,k}^i-F_i$ is connected. Since $k\geq 3$,
$A_{n,k}^i$ is not a complete graph, and so $A_{n,k}^i-F_i$ has at
least two vertices.
%Let $x$ and $y$ be two distinct vertices in $A_{n,k}^i-F_i$.
Since every vertex in $A_{n,k}^i-F_i$ has exactly $n-k$ outer
neighbors in $A_{n,k}-A_{n,k}^i$ and $|F-F_i|=n-k<2(n-k)$, at least
one of these outer neighbors is fault-free, and so $A_{n,k}-F$ is
still connected, a contradiction. Therefore, $A_{n,k}^i-F_i$ is
disconnected.

Let $H^i$ be a minimum component of $A_{n,k}^i-F_i$. Since $F$ is
a minimum separating set in $A_{n,k}$ and $F_i\subset F$, $H^i$
must be contained in some component $H$ in $A_{n,k}-F$. Note that
every vertex in $H^i$ has exactly $n-k$ outer neighbors in
$A_{n,k}-A_{n,k}^i$, each of them is in different $A_{n,k}^j$ with
$j\ne i$, and $A_{n,k}-A_{n,k}^i$ is $(n-k+1)$-connected. To
separate $H$ from $A_{n,k}-F$ by using $n-k$ vertices in $F-F_i$,
$H$ must be a single vertex, say $x$, and $F-F_i$ must be the
$(n-k)$ outer neighbors of $x$ in $A_{n,k}-A_{n,k}^i$. Thus,
$H=H^i=\{x\}$ and $F=N(x)$. Since $A_{n,k}-(A_{n,k}^i\cup
(F-F_i))$ is connected and every vertex in
$A_{n,k}^i-(F_i\cup\{x\})$ has $n-k$ fault-free outer neighbors in
$A_{n,k}-A_{n,k}^i$, $A_{n,k}-(F\cup\{x\})$ is connected.

Thus, when $n\geq 4$ and $k\geq 3$, $A_{n,k}$ is tightly super
$k(n-k)$-connected. The theorem follows.
\end{pf}

\vskip6pt

Since $A_{n,n-1}$ is isomorphic to a star graph $S_n$ and
$A_{n,n-2}$ is isomorphic to the alternating group graph $AG_n$,
by Theorem~\ref{thm3.1}, we have the following corollaries
immediately.

\vskip6pt

\begin{cor} \textnormal{(Cheng and Lipman~\cite{cl02})}
The star graph $S_n$ is tightly super $(n-1)$-connected for $n\geq
4$.
\end{cor}

\begin{cor}
The alternating group network $AG_n$ is tightly super
$(2n-4)$-connected for $n\geq 5$.
\end{cor}

\vskip6pt

In the rest of this section, we will investigate the fault tolerance
of $A_{n,k}$ when we remove a set $F$ of vertices, where $|F|$ is
roughly twice or three times of the traditional connectivity.

Let
 $$
 \begin{array}{c}
 I=\{i\in\langle n\rangle:\ f_i\geq (k-1)(n-k)\},\\
  A_{n,k}^I=\bigcup\limits_{i\in
 I}A_{n,k}^i, \ \  F_I=\bigcup\limits_{i\in I}F_i,
 \end{array}
 $$
and let
 $$
 J=\langle n\rangle\setminus I,\  A_{n,k}^J=\bigcup_{j\in J}A_{n,k}^j,\ F_J=\bigcup_{j\in J}F_j.
 $$

\begin{lem}\label{lem3.4}\
Let $F$ be a set of faulty vertices in $A_{n,k}$ with $|F|\leq
(3k-2)(n-k)-3$ and $k\geq 3$. Then $A^J_{n,k}-F_J$ is connected.
\end{lem}

\begin{pf}
If $|J|=0$ then there is nothing to do, and so assume $|J|\geq 1$.
By the hypothesis, for any $j\in J$, $f_j\leq (k-1)(n-k)-1$, that
is, $A^j_{n,k}-F_j$ is connected since $A^j_{n,k}$ is
$(k-1)(n-k)$-connected. Thus, if $|J|=1$ then the lemma holds.
Assume $|J|\geq 2$ below. To prove the lemma, we only need to show
that $A_{n,k}^i$ and $A_{n,k}^j$ are connected in $A_{n,k}^J-F_J$
for any two distinct $i,j\in J$. By (\ref{e2.1}), we have that
 \begin{equation}\label{e3.3}
 \begin{array}{rl}
 |E(i,j)| & =(n-2)(n-3)\cdots(n-k)\\
 &\left\{
 \begin{array}{ll}
 > 2((k-1)(n-k)-1)\ &\ {\rm if}\ k\ge 4\ {\rm or}\ k=3\ {\rm and}\ n\geq 6;\\
 =2(2n-7)\ &\ {\rm if}\ k=3\ {\rm and}\ n\in\{4,5\}.
 \end{array}\right.
 \end{array}
 \end{equation}

Thus, if there is a fault-free edge $e$ in $E(i,j)$, then
$A_{n,k}^i-F_i$ and $A_{n,k}^j-F_j$ can be connected by the
fault-free edge $e$ in $E(i,j)$. If there are no fault-free edges in
$E(i,j)$ then, by (\ref{e3.3}), $k=3$, $n\in\{4,5\}$ and
$f_i=f_j=2n-7$. In this case, $|F|=7n-24$, $|J|\geq 3$ and, for any
three distinct $i,j,x\in J$,
 \begin{equation}\label{e3.4}
 \begin{array}{rl}
 |F|-(f_i+f_j)\leq &(7n-24)-2(2n-7)\\
 =& 3n-10\\
 =&\left\{\begin{array}{ll}
  5<|E(i,x)|=|E(x,j)|\ &\ {\rm if}\ n=5;\\
  2=|E(i,x)|=|E(x,j)|\ &\ {\rm if}\ n=4.
 \end{array}\right.
 \end{array}
 \end{equation}

If $n=5$ then, by (\ref{e3.4}), there are a fault-free edge $e_1$ in
$E(i,x)$ and a fault-free edge $e_2$ in $E(x,j)$. Then $A_{5,3}^i$
and $A_{5,3}^j$ can be connected in $A_{5,3}-F$ by $A_{5,3}^x$ and
the fault-free edges $e_1$ and $e_2$.

If $n=4$, then $f_i=f_j=1$, and every vertex in $A_{4,3}^x$ has only
one outer neighbor for each $x\in\{1,2,3,4\}$. Thus, by
(\ref{e3.4}), there are a fault-free edge $e_1$ in $E(i,x)$ and a
fault-free edge $e_2$ in $E(x,j)$. Then $A_{4,3}^i$ and $A_{4,3}^j$
can be connected in $A_{4,3}-F$ by $A_{4,3}^x$ and the fault-free
edges $e_1$ and $e_2$.

 The lemma follows.\end{pf}

\vskip6pt

\begin{cor}\label{cor3.5}\
Let $F$ be a separating set of $A_{n,k}$ and $k\geq 3$. Then
  \begin{equation}\label{e3.5}
 1\leq |I|\leq \ \left\{\begin{array}{ll}
 2\ \ {\rm if}\ |F|\leq (2k-1)(n-k)-1;\\
 3\ \ {\rm if}\ |F|\leq (3k-2)(n-k)-3.
 \end{array}\right.
 \end{equation}
\end{cor}

\begin{cor}\label{cor3.6}\
Let $F$ be a separating set of $A_{n,k}$ with $|F|\leq
(3k-2)(n-k)-3$ and $k\geq 3$. If $H$ is a union of components of
$A_{n,k}-F$ that contain no vertices in $A^J_{n,k}-F$, then
 \begin{equation}\label{e3.6}
 N^{I}(H)\subseteq F_I\ \ {\rm and }\ N^{\overline{I}}(H)\subseteq F\setminus F_I.
 \end{equation}
\end{cor}

\vskip6pt

\begin{lem}\label{lem3.7}\
Let $F$ be a separating set of $A_{n,k}$ with $|F|\leq
(3k-2)(n-k)-3$ and $k\geq 3$. If there is some $i\in \langle
n\rangle$ such that $|F|-f_i\leq 2(n-k)-1$, then $A_{n,k}-F$ has
exactly two components, one of which is a single vertex.
\end{lem}

\begin{pf}
By the hypothesis, for any $j\in\langle n\rangle\setminus\{i\}$,
 $$
 f_j\leq |F|-f_i\leq 2(n-k)-1.
 $$
Since $|I|\geq 1$ by Corol1ary~\ref{cor3.5}, we have $I=\{i\}$.
Since $A_{n,k}-F$ is disconnected, and $A_{n,k}-(A^i_{n,k}\cup F)$
is connected by Lemma~\ref{lem3.4}, there is a component of
$A_{n,k}-F$ that contains no vertices in $A^J_{n,k}-F_J$. Let $H$ be
a union of such components of $A_{n,k}-F$. By
Corollary~\ref{cor3.6}, $N^{\overline{i}}(H)\subseteq F\setminus
F_i$. By (\ref{e2.2}) we have that
 $$%\begin{equation}\label{e3.6}
 |V(H)|(n-k)\leq |F|-f_i\leq 2(n-k)-1,
 $$%\end{equation}
which yields $|V(H)|\leq 1$, that is, $H$ is a single vertex, say
$u$. By the choice of $H$, other components of $A_{n,k}-F$ must
contain vertices in $A^J_{n,k}-F_J$. Since $A_{n,k}^J-F_J$ is
connected, $A_{n,k}-(F\cup\{u\})$ is connected. It follows that
$A_{n,k}-F$ has exactly two components, one of which is a single
vertex.

The lemma follows.
\end{pf}

\begin{lem}\label{lem3.8}\
Let $F$ be a separating set of $A_{n,k}$ with $|F|\leq
(3k-2)(n-k)-3$ and $k\geq 3$, and let $H$ be a subgraph of
$A_{n,k}^i-F_i$ for some $i\in \langle n\rangle$. If
$N_{A^i_{n,k}}(H)\subseteq F_i$, then $|V(H)|\leq 2$.
\end{lem}

\begin{pf}\
Let $h=|V(H)|$. We want to prove $h\leq 2$. Suppose to the contrary
that $h\geq 3$. Take a subset $T\subseteq V(H)$ with $|T|=3$. Let
$T'=V(H-T)$. By the hypothesis, $N_{A^i_{n,k}}(T)\setminus
T'\subseteq F_i$. Note that $A^i_{n,k}$ is $(k-1)(n-k)$-regular.

When $n=k+1$, by (\ref{e2.5}), any two vertices of $T$ have at most
one common neighbor in $A_{n,k}$. It follows that
 $$
 |N_{A^i_{n,k}}(T)|\geq 3(k-1)(n-k)-4.
 $$

When $n\geq k+2$, we denote $T=\{x,y,z\}$, and discuss as follows.

If $H[T]$ has no edges, then every pair of vertices in $T$ has at
most two common neighbors by (\ref{e2.5}), and so
 $$
 |N_{A^i_{n,k}}(T)|\geq 3(k-1)(n-k)-6.
 $$

If $H[T]$ has only one edge, say $e=(x,y)$, then $x$ and $y$ have
$n-k-1$ common neighbors, $z$ and $x$ (resp. $y$) have at most two
common neighbors by (\ref{e2.5}). It follows that
 $$
 |N_{A^i_{n,k}}(T)|\geq 3(k-1)(n-k)-(n-k-1)-6.
 $$

Similarly, by (\ref{e2.5}), we can obtain that if $H[T]$ has two
edges then
 $$
 |N_{A^i_{n,k}}(T)|\geq 3(k-1)(n-k)-2(n-k-2)-5;
 $$
if $H[T]$ has three edges,
 $$
 |N_{A^i_{n,k}}(T)|\geq 3(k-1)(n-k)-2(n-k-2)-6.
 $$

Summing all cases, we have that
\begin{eqnarray*}
 f_i&\geq & |N_{A^i_{n,k}}(T)\setminus T'|\\
      &\geq & |N_{A^i_{n,k}}(T)|-(h-3)\\
      &\geq & 3(k-1)(n-k)-6-2(n-k-2)-(h-3)\\
      &= & (3k-5)(n-k)-h+1,
\end{eqnarray*}
that is,
  \begin{equation}\label{e3.7}
  f_i\geq (3k-5)(n-k)-h+1.
 \end{equation}

Since $N^{\overline{i}}(H)\subseteq F-F_i$ by
Corollary~\ref{cor3.6}, $|F|-f_i\geq h(n-k)$, from which we have
that
 $$
 \begin{array}{rl}
 f_i\leq & |F|-h(n-k)\\
   \leq & (3k-2)(n-k)-3-h(n-k)\\
   = & (3k-2-h)(n-k)-3,
   \end{array}
   $$
that is,
  \begin{equation}\label{e3.8}
  f_i\leq (3k-2-h)(n-k)-3.
 \end{equation}
Combining (\ref{e3.7}) with (\ref{e3.8}), we have can deduce that
$(h-3)(n-k)\leq h-4$, a contradiction. Thus, we have $h\leq 2$. The
lemma follows.
\end{pf}

\vskip6pt

\begin{thm}\label{thm3.9}\
Let $F$ be a set of faulty vertices in $A_{n,k}$ with $|F|\leq
(2k-1)(n-k)-1$ and $k\geq 3$. If $A_{n,k}-F$ is disconnected, then
it has exactly two components, one of which is a single vertex or a
single edge.
\end{thm}

\begin{pf}\
Since $A_{n,k}-F$ is disconnected, $F$ is a separating set of
$A_{n,k}$.
% and $|F|\geq\kappa(A_{n,k})=k(n-k)$. If $|F|=k(n-k)$
%then, by Theorem~\ref{thm3.1}, $A_{n,k}-F$ has exactly two
%components, one of which is a single vertex, and so the theorem
%holds. Thus, in the following discussion, we assume that
% $$%\begin{equation}\label{e3.8}
% k(n-k)+1\leq |F|\leq (2k-1)(n-k)-1.
% $$%\end{equation}
%Under this hypothesis, we want to prove that $A_{n,k}-F$ has
%exactly two components, one of which is either a single vertex or
%a single edge.

Suppose that there exists some $i\in\langle n\rangle$ such that
$f_i\geq (2k-3)(n-k)$. Since $|F|\leq (2k-1)(n-k)-1\leq
(3k-2)(n-k)-3$ and
 $$
|F|-f_i\leq ((2k-1)(n-k)-1)-(2k-3)(n-k)=2(n-k)-1,
 $$
$A_{n,k}-F$ has exactly two components, one of which is a single
vertex by Lemma~\ref{lem3.7}.

We now assume that $f_i\leq (2k-3)(n-k)-1$ for any $i\in\langle
n\rangle$. Then
 $$
 \begin{array}{rl}
  |V(A^i_{n,k}-F_i)|=& (n-1)(n-2)\cdots  (n-k)-f_i\\
 \geq & (n-1)(n-2)\cdots (n-k)-((2k-3)(n-k)-1)\\
 \geq &2.
 \end{array}
 $$

Since
 $$
 |F|\leq (2k-1)(n-k)-1\leq
 (3k-2)(n-k)-3,
 $$
by Lemma~\ref{lem3.4} $A_{n,k}^J-F_J$ is connected. Let $H$ be a
union of components of $A_{n,k}-F$ that contain no vertices in
$A^J_{n,k}-F_J$. Thus, $H$ is in $A_{n,k}^I$. By the choice of
$H$, other components of $A_{n,k}-F$ must contain vertices in
$A^J_{n,k}-F_J$. Since $A_{n,k}^J-F_J$ is connected,
$A_{n,k}-(F\cup V(H)$ is connected. Thus, to complete the proof of
the theorem, we only need to show that $H$ is either a single
vertex or a single edge. Consider two cases according to $|I|=1$
or $|I|=2$ by Corollary~\ref{cor3.5}.

\vskip6pt

{\bf Case 1.}\ $|I|=1$, and let $I=\{i\}$.

Let $h=|V(H)|$. Then $h\leq 2$ by Lemma~\ref{lem3.8}. If $h=1$, then
$H$ is a single vertex.

If $h=2$, we want to prove that $H$ is a single edge. Suppose to
the contrary that $H$ consists of two isolated vertices, say $u$
and $v$. Then $u$ and $v$ are not adjacent, $N(u)\cup
N(v)\subseteq F$. By (\ref{e2.5}), we deduce a contradiction as
follows.
 $$
 \begin{array}{rl}
 |F|\geq & |N(u)\cup N(v)|=|N(u)|+|N(v)|-|N(u)\cap N(v)|\\
 =& 2k(n-k)-|N(u)\cap N(v)|\\
 >&(2k-1)(n-k)-1\geq |F|.
 \end{array}
 $$
Thus, $H$ is a single edge.

\vskip6pt

{\bf Case 2.}\ $|I|=2$, and let $I=\{i,j\}$.

Under our hypothesis, by (\ref{e2.2}) and (\ref{e3.6}), we have that
 $$%\begin{equation}\label{e3.9}
 \begin{array}{rl}
 n-k-1\leq | N^{\overline{I}}(H)|\leq &|F\setminus(F_i\cup F_j)|\\
 \leq &(2k-1)(n-k)-1-2((k-1)(n-k))\\
 =&n-k-1.
 \end{array}
 $$% \end{equation}
Thus, $| N^{\overline{I}}(H)|=n-k-1$.

Thus, by (\ref{e2.2}) and (\ref{e2.2}), there is exactly one vertex
in $(A_{n,k}^i-F_i)\cap V(H)$ such that exact one of its outer
neighbors is in $A_{n,k}^j$ and others are in $F\setminus(F_i\cup
F_j)$. Similarly, there is exactly one vertex in
$(A_{n,k}^j-F_j)\cap V(H)$ such that exact one of its outer
neighbors is in $A_{n,k}^i$ and others are in $F\setminus(F_i\cup
F_j)$. Thus, $H$ is a single edge.

\vskip6pt

%Summing our discussion above, whichever $|I|=1$ or $|I|=2$, we have
%$A_{n,k}-F$ have exactly two components, one of them is an isolated
%vertex or a single edge.
The proof of the theorem is complete.
\end{pf}

Since $A_{n,n-1}$ is isomorphic to a star graph $S_n$ and
$A_{n,n-2}$ is isomorphic to a alternating group graph $AG_n$, by
Theorem~\ref{thm3.9}, we have the following corollaries
immediately.

\begin{cor} \label{cor3.10}\textnormal{(Cheng and Lipman~\cite{cl02})}
Let $F$ be a set of faulty vertices in the star graph $S_n$ with
$|F|\leq 2n-4$ and $n\geq 4$. If $S_n-F$ is disconnected, then it
has exactly two components, one of which is either a single
vertex, or a single edge.
\end{cor}

\begin{cor}\label{cor3.11}
Let $F$ be a set of faulty vertices in the alternating group graph
$AG_n$ with $|F|\leq 4n-11$ and $n\geq 5$. If $AG_n-F$ is
disconnected, then it has exactly two components, one of which is
either a single vertex, or a single edge.
\end{cor}

We now discuss the fault tolerance of $A_{n,k}$ with more faulty
vertices up to $(3k-2)(n-k)-4$ when $n\geq k+2$ and
$(3k-2)(n-k)-3$ when $n=k+1$, where the latter we write $3n-8$ for
$(3k-2)(n-k)-3$.

\begin{thm}\label{thm3.12}\
Let $F$ be a set of faulty vertices in $A_{n,k}$ ($k\geq 4$) with
$|F|\leq (3k-2)(n-k)-4$ when $n\geq k+2$ and $|F|\leq 3n-8$ when
$n=k+1$. If $A_{n,k}-F$ is disconnected, then it either has two
components, one of which is an isolated vertex or an isolated edge,
or has three components, two of which are isolated vertices.
\end{thm}

\begin{pf}\
Since $A_{n,k}-F$ is disconnected, $F$ is a separating set of
$A_{n,k}$.
% and $|F|\geq\kappa(A_{n,k})=k(n-k)$. If $|F|\leq
%(2k-1)(n-k)-1$, then the theorem holds by Theorem~\ref{thm3.9}. So,
%in the following discussion, we assume that
% $$
%(2k-1)(n-k)\leq |F|\leq (3k-2)(n-k)-3.
% $$

If there exists some $i\in\langle n\rangle$ such that
 $$
 |F|-f_i \leq 2(n-k)-1,
 $$
by Lemma~\ref{lem3.7}, $A_{n,k}-F$ has exactly two components, one
of which is a single vertex, and so the theorem holds. We now assume
that, for any $i\in\langle n\rangle$
 $$
 f_i\leq
 \left\{\begin{array}{ll}
 (3k-4)(n-k)-4\ \ & {\rm for}\  n\geq k+2;\\
 (3k-4)(n-k)-3\ \ & {\rm for}\  n=k+1.
 \end{array}\right.
 $$
Then $|V(A^i_{n,k}-F_i)|\geq 2$.

Let $H$ be a union of components of $A_{n,k}-F$ that contain no
vertices in $A^J_{n,k}-F_J$, and let $h=|V(H)|$. By
Lemma~\ref{lem3.4}, $A_{n,k}^J-F_J$ is connected. Thus, $H$ is in
$A_{n,k}^I$. By the choice of $H$, other components of $A_{n,k}-F$
must contain vertices in $A^J_{n,k}-F_J$. Since $A_{n,k}^J-F_J$ is
connected, $A_{n,k}-(F\cup V(H))$ is connected. Thus, to complete
the proof of the theorem, we only need to show that $h\leq 2$.

By Corollary~\ref{cor3.5}, $1\leq |I|\leq 3$. If $|I|=3$, under our
hypothesis, we have that
 \begin{equation}\label{e3.9}
 |F\setminus F_I|\leq (3k-2)(n-k)-4-3((k-1)(n-k))=n-k-4.
 \end{equation}
Thus, by (\ref{e2.2}), (\ref{e3.3}) and (\ref{e3.9}), we can deduce
a contradiction as follows.
 $$
 n-k-2\leq |F\setminus F_I|\leq n-k-4.
 $$
Thus,  $1\leq |I|\leq 2$. If $|I|=1$, then $h\leq 2$ by
Lemma~\ref{lem3.8}. We only need to consider the case of $|I|=2$.
Let $I=\{i,j\}$, and let $h_i$ and $h_j$ be the numbers of vertices
of $H$ that lie in $A_{n,k}^i$ and  $A_{n,k}^j$, respectively. Then
$h_i\leq 2$ and $h_j\leq 2$ by Lemma~\ref{lem3.8}. Without loss of
generality, assume $h_i\geq h_j$ and $f_i\geq f_j$.

Note that $A_{n,k}^i$ is isomorphic to $A_{n-1,k-1}$. If $f_i\leq
(2k-3)(n-k)-2$ then, when $k\geq 4$, applying Theorem~\ref{thm3.9}
to $A_{n,k}^i$, we have $h_i\leq 1$ since $F_i$ can not isolate an
edge from $A_{n,k}^i$. Thus, $h=h_i+h_j\leq 2$. So, in the following
discussion, we assume that
 \begin{equation}\label{e3.10}
 f_i\geq (2k-3)(n-k)-1.
 \end{equation}

If $f_j\geq k(n-k)$, when $|F|\leq (3k-2)(n-k)-3$, we have that
 $$
 \begin{array}{rl}
 |F\setminus F_I|\leq & (3k-2)(n-k)-3-k(n-k)-(2k-3)(n-k)+1\\
 = &n-k-2.
 \end{array}
 $$
Note that, for every vertex of $V(H)\cap V(A_{n,k}^i)$, it has at
most one outer neighbor in $A_{n,k}^j$ and others in $F\setminus
F_I$. By (\ref{e2.2}) and (\ref{e3.3}), we have that
 $$
 h_i(n-k-1)\leq |F\setminus F_I|\leq n-k-2,
 $$
which implies $h_i=0$, and so $h=h_i+h_j\le 2$.

Thus, under the condition (\ref{e3.10}), the remainder of the proof
is to consider the case that
 \begin{equation}\label{e3.11}
 (k-1)(n-k)\leq f_j\leq k(n-k)-1.
 \end{equation}

We first note that, when $f_j\leq k(n-k)-1$,
 $$
 f_j\leq k(n-k)-1\leq (2k-3)(n-k)-2.
 $$
Thus, $F_j$ isolates at most one vertex in $A_{n,k}^j$ by
Theorem~\ref{thm3.9}, that is, $h_j\leq 1$. If $h_j=0$, then $h\leq
2$, and so the theorem holds. Assume $h_j=1$ below.

By  the condition (\ref{e3.10}) and  the condition (\ref{e3.11}), we
have that
\begin{equation}\label{e3.12}
\begin{array}{rl}
 f_i+f_j &\geq (2k-3)(n-k)-1+(k-1)(n-k)\\
         & =(3k-4)(n-k)-1.
         \end{array}
 \end{equation}
Thus, by (\ref{e3.12}), when $|F|\leq (3k-2)(n-k)-4$ and $n\geq
k+2$,
\begin{equation}\label{e3.13}
 |F\setminus F_I|=|F|-f_i-f_j\leq 2(n-k)-3,
 \end{equation}
and when $|F|\leq 3n-8$ and $n=k+1$,
\begin{equation}\label{e3.14}
 |F\setminus F_I|=|F|-f_i-f_j=0.
 \end{equation}

Suppose to the contrary that $h_i=2$. Let
 $$
 V(H)\cap V(A_{n,k}^i)=\{x,y\}\ \ {\rm and}\ \ V(H)\cap V(A_{n,k}^j)=\{z\}.
 $$
Then at least one of $x$ and $y$ is not adjacent to $z$ by
(\ref{e2.2}). Without loss of generality, let $x$ be not adjacent to
$z$. Then, by (\ref{e2.3}),
\begin{equation}\label{e3.15}
 |N^{\overline{I}}(x)\cap N^{\overline{I}}(z)|=0.
 \end{equation}

By (\ref{e2.2}), we have
\begin{equation}\label{e3.16}
 |N^{\overline{I}}(x)\cap N^{\overline{I}}(y)|=0,
 \end{equation}

and

\begin{equation}\label{e3.17}
 |N^{\overline{I}}(x)\cap N^{\overline{I}}(y)\cap N^{\overline{I}}(z)|=0,
 \end{equation}

 When $n\geq k+2$, considering outer neighbors of $y$ and $z$,
by (\ref{e2.3}) and (\ref{e2.5}), we have that
 \begin{equation}\label{e3.18}
 |N^{\overline{I}}(y)\cap N^{\overline{I}}(z)|=
 \left\{\begin{array}{ll}
  0 \ \ & {\rm if}\ (y,z)\notin E(A_{n,k});\\
  n-k-1\ \ & {\rm if}\ (y,z)\in E(A_{n,k}).
 \end{array}\right.
 \end{equation}
By (\ref{e3.6}), (\ref{e3.15})-(\ref{e3.18}), we have that
 $$
 \begin{array}{rl}
 |F\setminus F_I| &\geq |N^{\overline{I}}(H)|\\
  &=\sum\limits_{u\in \{x,y,z\}}|N^{\overline{I}}(u)|-
  \sum\limits_{u\neq v\in \{x,y,z\}}|N^{\overline{I}}(u)\cap N^{\overline{I}}(v)|\\
  &\qquad \qquad \qquad \qquad -|N^{\overline{I}}(x)\cap N^{\overline{I}}(y)\cap N^{\overline{I}}(z)|\\
  &\geq \left\{\begin{array}{ll}
 2(n-k-1)\ \ & {\rm if}\ (y,z)\in E(A_{n,k});\\
 3(n-k-1)\ \ & {\rm if}\ (y,z)\notin  E(A_{n,k}),
 \end{array}\right.
 \end{array}
 $$
which contradicts (\ref{e3.13}). Thus, $h_i\leq 1$ and so
$h=h_i+h_j\le 2$.

When $n=k+1$, (\ref{e3.14}) implies that $|F|=3n-8$, $f_j=n-2$ and
$f_i=2n-6$. In other words, $F_j=N_{A_{n,k}^j}(z)$ and
$F_i=N_{A_{n,k}^i}(x,y)$, the latter implies that $x$ and $y$ are
adjacent. Since $n=k+1$, the only outer neighbor of $x$, say $u$,
and the only outer neighbor of $y$, say $v$, must be in
$F_j\cup\{z\}$. Similarly, the only outer neighbor of $z$, say
$w$, must be in $F_i\cup\{x,y\}$. Since $x$ is not adjacent to
$z$, $u\in F_j$. If $v=z$ then $|N(x)\cap N(z)|=2$, which
contradict (\ref{e2.5}). Assume $v\in F_j$ below. If $w\in N(x)$,
then $|N(x)\cap N(z)|=2$; if $w\in N(y)$, then $|N(y)\cap
N(z)|=2$. No matter which case, it contradicts (\ref{e2.5}).

The proof of the theorem is complete.
\end{pf}

The theorem~\ref{thm3.12} is optimal in the following sense. When
$n\geq k+2$ and $k\geq 4$, we select such three vertices $x,y\in
V(A_{n,k}^i)$ and $z\in V(A_{n,k}^j)$ that $(y,z)\in E(i,j)$ and
$(x,y)\in E(A^i_{n,k})$, see Figure~\ref{f2}. Set $F=N(x,y,z)$. By
(\ref{e2.5}), we have that
 $$
 \begin{array}{rl}
 |N(x)\cap N(y)|=|N(y)\cap N(z)|=n-k-1, |N(x)\cap N(z)|=2.
 \end{array}
 $$
Then
 $$
 |F|=3k(n-k)-2(n-k-1)-5=(3k-2)(n-k)-3.
 $$
$A_{n,k}-F$ is connected and contains a path of length three.

\begin{center}
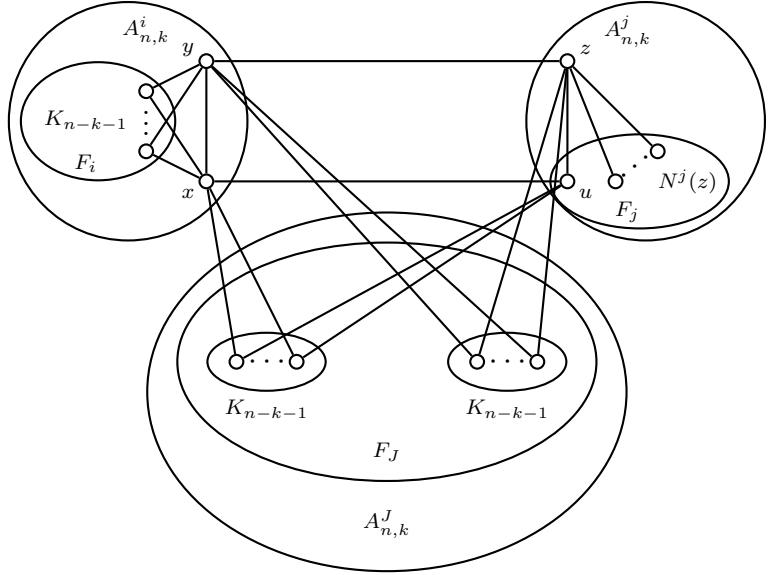
\begin{figure}[ht]
\psset{unit=.8}
\begin{pspicture}(-10,-4)(8,6)

\cnode(-1.5,0){3pt}{1}  \cnode(-2.5,0){3pt}{2} \rput(-2,0){$\cdots$}
\psellipse(-2,0)(1,.5)  \rput(-2,-.8){\scriptsize $K_{n-k-1}$}
\cnode(1.5,0){3pt}{3}  \cnode(2.5,0){3pt}{4} \rput(2,0){$\cdots$}
\psellipse(2,0)(1,.5)  \rput(2,-.8){\scriptsize $K_{n-k-1}$}
\psellipse(0,0)(3.5,2)  \rput(0,-1.5){\scriptsize $F_J$}

\cnode(-3,3){3pt}{5}\rput(-3.3,2.8){\scriptsize $x$}
\cnode(-3,5){3pt}{6}\rput(-3.3,5.2){\scriptsize $y$} \ncline{5}{6}
\cnode(3,3){3pt}{7}\rput(3.3,2.8){\scriptsize $u$}
\cnode(3,5){3pt}{8}\rput(3.3,5.2){\scriptsize $z$} \ncline{7}{8}
\ncline{5}{7} \ncline{6}{8}

\ncline{1}{5} \ncline{1}{7}  \ncline{2}{5}  \ncline{2}{7}
\ncline{3}{6} \ncline{3}{8}  \ncline{4}{6}  \ncline{4}{8}

\cnode(-4,3.5){3pt}{9} \cnode(-4,4.5){3pt}{10}
\rput(-4,4.1){$\vdots$} \rput(-5,4){\scriptsize $K_{n-k-1}$}
\rput(-5,3.3){\scriptsize $F_i$}   \psellipse(-4.8,4)(1.3,1)
\ncline{5}{9} \ncline{6}{9}  \ncline{5}{10}  \ncline{6}{10}

\cnode(3.8,3){3pt}{11} \cnode(4.5,3.5){3pt}{12}
\rput[b]{40}(4.2,3.15){$\cdots$} \rput(5,3){\scriptsize $N^j(z)$}
\rput(4,2.5){\scriptsize $F_j$}   \psellipse(4.2,3)(1.5,.8)
\ncline{8}{11} \ncline{8}{12}

\rput(-4,5.5){\scriptsize $A_{n,k}^i$}   \psellipse(-4.3,4)(2,2)
\rput(4,5.5){\scriptsize $A_{n,k}^j$}   \psellipse(4.3,4)(2,2)
\rput(0,-2.7){\scriptsize $A_{n,k}^J$}   \psellipse(0,-.5)(4,3)

\end{pspicture}
\caption{\label{f2}\footnotesize{The distribution of fault set $F$
in $A_{n,k}$ with $k\geq 4$ and $n\geq k+2$}}
\end{figure}
\end{center}

Since $A_{n,n-1}$ is isomorphic to a star graph $S_n$ and
$A_{n,n-2}$ is isomorphic to a alternating group graph $AG_n$, by
Theorem~\ref{thm3.12}, we have the following corollaries
immediately.

\begin{cor}\label{cor3.13} \textnormal{(Cheng and Lipt\'ak~\cite{cl06})}
Let $F$ be a set of faulty vertices in the star graph $S_n$ with
$|F|\leq 3n-8$ and $n\geq 5$. If $S_n-F$ is disconnected, then it
it either has two components, one of which is an isolated vertex
or an edge, or has three components, two of which are isolated
vertices.
\end{cor}

\begin{cor} \label{cor3.14}
Let $F$ be a set of faulty vertices in the alternating group
network $AG_n$ with $|F|\leq 6n-20$ and $n\geq 6$. If $AG_n-F$ is
disconnected, then it either has two components, one of which is
an isolated vertex or an edge, or has three components, two of
which are isolated vertices.
\end{cor}

\section{\it\bf Diagnosability of arrangement graph}

The comparison diagnosis strategy of a graph $G=(V,E)$ can be
modeled as a multi-graph $M=(V,C)$, where $C$ is a set of labelled
edges. If the processors $u$ and $v$ can be compared by the
processor $w$, there exists an labelled edge $(u,v)$ in $C$,
denoted by $(u,v)_w$. We call $w$ the comparator of $u$ and $v$.
Since different comparators can compare the same pair of
processors, $M$ is a multi-graph. Denote the comparison result as
$\sigma ((u,v)_w)$ such that $\sigma ((u,v)_w)=0$ if the outputs
of $u$ and $v$ agree, and $\sigma ((u,v)_w)=1$ if the outputs
disagree. If the comparator $w$ is fault-free and $\sigma
((u,v)_w)=0$, the processors $u$ and $v$ are fault-free; while
$\sigma ((u,v)_w)=1$, at least one of the three processors $u$,
$v$ and $w$ is faulty. The collection of the comparison results
defined as a function $\sigma :\ C\rightarrow \{0,1\}$, is called
the {\it syndrome} of the diagnosis. If the comparator $w$ is
faulty, the comparison result is unreliable. A faulty comparator
can lead to unreliable results, so a set of faulty vertices may
produce different syndromes. A subset $F\subsetneq V$ is said to
be {\it compatible} with a syndrome $\sigma$ if $\sigma$ can arise
from the circumstance that all vertices in $F$ are faulty and all
vertices in $V-F$ are fault-free.  A system $G$ is said to be {\it
diagnosable} if, for every syndrome $\sigma$, there is a unique
$F\subset V$ that is compatible with $\sigma$. A system is said to
be a $t$-diagnosable if the system is diagnosable as long as the
number of faulty vertices does not exceed $t$. The maximum number
of faulty vertices that the system $G$ can guarantee to identify
is called the {\it diagnosability} of $G$, write as $t(G)$. Let
$\sigma_F=\{\sigma \ |\ \sigma$ is compatible with $F\}$. Two
distinct subsets $F_1$ and $F_2$ of $V(G)$ are said to be {\it
indistinguishable} if and only if $\sigma_{F_1}\cap
\sigma_{F_1}\neq\phi$, and distinguishable
otherwise~\cite{hcsht09, s92, lthcl08}. There are several
different ways to verify whether a system is $t$-diagnosable under
the comparison approach. The following lemma obtained by Sengupta
and Dahbura~\cite{s92} gives necessary and sufficient conditions
to ensure distinguishability.

\begin{lem}\label{lem4.1} \textnormal{(Sengupta and Dahbura
~\cite{s92})}\ Let $G$ be a graph, $F_1$ and $F_2$ be two distinct
subsets of vertices in $G$. The pair $(F_1,F_2)$ is distinguishable
if and only if at least one of the following conditions is
satisfied.

(1) There are two distinct vertices $u$ and $w\in V(G-F_1\cup F_2)$
and a vertex $v\in F_1\Delta F_2$ such that $(u,v)_w\in C$, where
$F_1\bigtriangleup F_2=(F_1\setminus F_2)\cup (F_2\setminus F_1)$;

(2) There are two distinct vertices $u$ and $v\in F_1\setminus F_2$
(or $F_2\setminus F_1$) and a vertex $w\in V(G-F_1\cup F_2)$ such
that $(u,v)_w\in C$.
%;
%(3) There are two distinct vertices $u$ and $v\in F_2\setminus F_1$
%and a vertex $w\in V(G-F_1\cup F_2)$ such that $(u,v)_w\in C$.
\end{lem}

Lin et al.~\cite{lthcl08} introduced the so-called conditional
diagnosability of a system under the situation that no set of
faulty vertices can contain all neighbors of any vertex in the
system. A fault-set $F\subset V(G)$ is called a conditional
fault-set if $G-F$ has no isolated vertex. A system $G(V,E)$ is
said to be conditionally $t$-diagnosable if $F_1$ and $F_2$ are
distinguishable for each pair ($F_1,F_2$) of distinct conditional
fault-sets in $G$ with $|F_1|\leq t$ and $|F_2|\leq t$. The {\it
conditional diagnosability} of $G$, denoted by $t_c(G)$ is defined
as the maximum value of $t$ for which $G$ is conditionally
$t$-diagnosable. Clearly, $t_c(G)\geq t(G)$. Zhou and
Xiao~\cite{zx10} obtained the conditional diagnosability of the
alternating group networks based on the fault tolerance of this
network structure. This section will focus on the conditional
diagnosability of arrangement graphs.

\begin{thm}\label{thm4.2}\
 $t_c(A_{n,k})\leq (3k-2)(n-k)-3$ for $k\geq 4$, $n\geq k+2$; $t_c(A_{n,k})\leq 3n-7$ for $k\geq 4$, $n=k+1$.
\end{thm}

\begin{pf}
When $n\geq k+2$, we select four vertices $x,y,z,u\in V(A_{n,k})$,
such that $(x,u), (y,z)\in E(i,j)$, and $(x,y)\in E(A^i_{n,k})$,
then $(u,z)\in E(A^j_{n,k})$. Set $A=N[x,y,z]$, $F_1=A-\{y,z\}$,
and $F_2=A-\{x,y\}$. We get
$$
|F_1|=|F_2|=(3k-2)(n-k)-2, {\rm and}\  |F_1-F_2|=|F_2-F_1|=1.
$$

It is easy to check that $F_1$ and $F_2$ are two conditional fault
sets, and $F_1$ and $F_2$ are indistinguishable. Thus, we have
$$t_c(A_{n,k})\leq (3k-2)(n-k)-3.$$

When $n=k+1$, we select three vertices $x,y,z\in V(A_{n,k})$, such
that $(x,y), (y,z)\in E(A_{n,k})$. By (~\ref{e2.5}), any two of
$x,y,z$ have no common neighbor. Set
$$
A=N[x,y,z], F_1=A-\{y,z\}, {\rm and} \ F_2=A-\{y,z\}.
$$

We get $|F_1|=|F_2|=3n-6$, and $|F_1-F_2|=|F_2-F_1|=1$. It is easy
to check that $F_1$ and $F_2$ are two conditional fault sets, and
$F_1$ and $F_2$ are indistinguishable. Thus, we have
$t_c(A_{n,k})\leq 3n-7$.

\end{pf}

\begin{lem}\label{lem4.3} \
 Let $F_1$ and $F_2$ be any two distinct conditional fault-sets
of $A_{n,k}$ with $|F_1|\leq (3k-2)(n-k)-3$, $|F_2|\leq
(3k-2)(n-k)-3$ for  $k\geq 4$, $n\geq k+2$; or $|F_1|\leq 3n-7$,
$|F_2|\leq 3n-7$ for  $k\geq 4$, $n=k+1$. Denote by $H$ the
maximum component of $A_{n,k}-F_1\cap F_2$. Then, for every vertex
$u\in F_1\Delta F_2$, $u\in H$.
\end{lem}

\begin{pf}\
Without loss of generality, we assume that $u\in F_1- F_2$. Since
$F_2$ is a conditional faulty set, there is a vertex
$v\in(A_{n,k}-F_2)-\{u\}$ such that $(u,v)\in E(A_{n,k})$. Suppose
that $u$  is not a vertex of $H$. Then $v$ is not in $H$, so $u$
and $v$ are in one small component of $A_{n,k}-F_1\cap F_2$. Since
$F_1$ and $F_2$ are distinct, we have
$$
|F_1\cap F_2|\leq (3k-2)(n-k)-4\  {\rm for}\ n\geq k+2;$$

or
$$
|F_1\cap F_2|\leq 3n-8\  {\rm for}\ n=k+1.\qquad\qquad\qquad$$

Hence $\{u,v\}$ forms a component $K_2$ in $A_{n,k}-F_1\cap F_2$
by Theorem~\ref{thm3.12}, i.e., the vertex $u$ is the unique
neighbor of $v$ in $A_{n,k}-F_1\cap F_2$. This is a contradiction
since $F_1$ is a conditional fault set, but all the neighbors of
$v$ are faulty in $A_{n,k}-F_1$.
\end{pf}

\begin{lem}\label{lem4.4} \textnormal{(C. K. Lin ~\cite{lthcl08})}
Let $G$ be a graph with $\delta (G)\geq 2$, and let $F_1$ and
$F_2$ be any two distinct conditional fault-sets of $G$ with
$F_1\subset F_2$. Then, $(F_1,F_2)$ is a distinguishable
conditional pair under the comparison diagnosis model.
\end{lem}

\begin{lem}\label{lem4.5} \
Let $F_1$ and $F_2$ be any two distinct conditional fault-sets of
$A_{n,k}$. If $|F_1|=(3k-2)(n-k)-3$ and $|F_2|=(3k-2)(n-k)-3$
$k\geq 4$, $n\geq k+2$; or $|F_1|\leq 3n-7$, $|F_2|\leq 3n-7$ for
$k\geq 4$, $n=k+1$. Then, $(F_1,F_2)$ is a distinguishable
conditional pair under the comparison diagnosis model.
\end{lem}

\begin{pf}\
By Lemma~\ref{lem4.4}, $(F_1,F_2)$ is a distinguishable
conditional pair if $F_1\subset F_2$ or $F_2\subset F_1$. Now, we
assume that $|F_1-F_2|\geq 1$, and $|F_2-F_1|\geq 1$. Let
$S=F_1\cap F_2$. Then we have $|S|\leq (3k-2)(n-k)-4$ for $k\geq
4$, $n\geq k+2$; or $|S|\leq 3n-8$ for  $k\geq 4$, $n=k+1$. Let
$H$ be the largest connected component of $A_{n,k}-F_1\cup F_2$.
By Lemma~\ref{lem4.3}, every vertex in $F_1\Delta F_2$ is in $H$.
\par We claim that $H$  has a vertex $u$ outside $F_1\cup F_2$ that has no neighbor in $H$.
Since every vertex has degree $k(n-k)$, the vertices in $S$ can have
at most $k(n-k)|S|$ neighbors in $H$. There are at most
$|F_1|+|F_2|-|S|$ vertices in $F_1\cup F_2$ and at most two vertices
of $A_{n,k}-S$ may not belong to $H$ by Theorem~\ref{thm3.12}. Thus,
we have:
$$
 \begin{array}{rl}
\frac{n!}{(n-k)!}&-k(n-k)|S|-(|F_1|+|F_2|-|S|)-2\\
&\geq \frac{n!}{(n-k)!}-(k(n-k)+1)\times ((3k-2)(n-k)-4)-4\\
         &\geq 4 \ {\rm \ for}\  k\geq 4, n\geq k+2;
\end{array}
$$

and
$$
 \begin{array}{rl}
\frac{n!}{(n-k)!}-& k(n-k)|S|-(|F_1|+|F_2|-|S|)-2 \qquad\qquad\quad\\
\geq & n!-n\times (3n-8)-2\\
\geq & n!-3n^2+8n-2\\
         \geq &4\ {\rm \ for}\  k\geq 4, n=k+1.
\end{array}
$$

Thus, there must be some vertex of $H$ outside $F_1\cup F_2$,
which has no neighbors in $S$. Let $u$ be such a vertex.
\par If $u$ has no neighbor in $F_1\cup F_2$, then we can find a path of length at least two
within $H$ to a vertex $v$ in $F_1\cup F_2$. We may assume that
$v$ is the first vertex of $F_1\Delta F_2$ on this path, and let
$q$ and $w$ be the two vertices on this path immediately before
$v$ (we may have $u=q$), so $q$ and $w$ are not in $F_1\cup F_2$.
The existence of the edges $(q,w)$ and $(w,v)$ ensures that
$(F_1,F_2)$ is a distinguishable conditional pair of $A_{n,k}$ by
Lemma~\ref{lem4.1}. Now we assume that $u$ has a neighbor in
$F_1\Delta F_2$. Since the degree of $u$ is at least 3, and $u$
has no neighbor in $S$, there are three possibilities:

(1) $u$ has two neighbors in $F_1\setminus F_2$; or

(2) $u$ has two neighbors in $F_2\setminus F_1$; or

(3) $u$ has at least one neighbor outside $F_2\cup F_1$.

In each sub-case above, Lemma~\ref{lem4.1} implies that
$(F_1,F_2)$ is a distinguishable conditional pair of $A_{n,k}$
under the comparison diagnosis model, and so the proof is
complete.
\end{pf}

Theorem~\ref{thm4.2} tells us that $t_c(A_{n,k})\leq
(3k-2)(n-k)-3$ for $k\geq 4$, $n\geq k+2$; $t_c(A_{n,k})\leq 3n-7$
for $k\geq 4$, $n=k+1$. Lemma~\ref{lem4.5} shows that
$t_c(A_{n,k})\geq (3k-2)(n-k)-3$ for $k\geq 4$, $n\geq k+2$;
$t_c(A_{n,k})\geq 3n-7$ for $k\geq 4$, $n=k+1$. Thus, we have the
following results.

\begin{thm}\label{thm4.6}\
$t_c(A_{n,k})=(3k-2)(n-k)-3$ for $k\geq 4$, $n\geq k+2$;
$t_c(A_{n,k})=3n-7$ for $k\geq 4$, $n=k+1$.
\end{thm}

Since $A_{n,n-1}$ is isomorphic to a star graph $S_n$ and
$A_{n,n-2}$ is isomorphic to a alternating group graph $AG_n$, by
Theorem~\ref{thm3.12}, we have the following corollaries
immediately.

\begin{cor} \textnormal{(C. K. Lin, et al.~\cite{lthcl08})}
The conditional diagnosability of the star graph $S_n$ under the
comparison model is $t_c(S_n)=3n-7$ for $n\geq 5$.
\end{cor}

\begin{cor}
The conditional diagnosability of the alternating group graph
$AG_n$ under the comparison model is $t_c(AG_n)=6n-19$ for $n\geq
6$.
\end{cor}

\section{\it\bf  Conclusion }
The paper derives the fault resiliency of arrangement graphs, and
then uses the fault resiliency to evaluate fault diagnosability of
the arrangement graphs under the comparison model. The fault
resiliency of the arrangement graphs may also reveal its conditional
connectivity of high order. This method can be also applied to other
complex network structure, such as $(n,k)$-star graphs.

{}}
\end{document}